\documentclass[12pt,a4paper]{amsart}
\usepackage{latexsym}
\usepackage{amsmath,amsthm,amssymb,amsfonts}

\allowdisplaybreaks

\newcommand{\head}{
{\noindent \small V International Meeting on Approximation Theory
of the University of Ja\'{e}n.\\ \'{U}beda, 9--14 June 2004}
\vspace{10pt}}

\newtheorem{theorem}{Theorem}

\begin{document}

\newcommand{\ZZ}{\mathbb{Z}}
\newcommand{\PP}{\mathbb{P}}
\newcommand{\RR}{\mathbb{R}}
\newcommand{\NN}{\mathbb{N}}
\newcommand{\abs}[1]{\left| #1 \right|}

\newtheorem{defin}[theorem]{Defin}
\head \vspace{1.5cm}

\title{RECENT RESULTS ON NEAR-BEST SPLINE QUASI-INTERPOLANTS}

\author{\underline{Paul Sablonni\`ere}}

\thanks {}

\maketitle

\markboth{P. Sablonni\`ere}{Near-best spline quasi-interpolants}

%*********************************************************  ABSTRACT
*************************************

\begin{center} \textbf{Abstract}\end{center}

Roughly speaking, a near-best (abbr. NB)  quasi-interpolant (abbr. QI) is
an approximation
operator of the form $Q_af=\sum_{\alpha\in A} \Lambda_\alpha (f) B_\alpha$
where the $B_\alpha$'s
are B-splines and the $\Lambda_\alpha (f)$'s are linear discrete or
integral forms acting on the
given function $f$. These forms  depend on a finite number of coefficients
which are the
components of vectors $a_\alpha$ for
$\alpha\in A$. The index $a$ refers to this sequence of vectors. In order
that $Q_a p=p$ for all
polynomials $p$ belonging to some subspace included in the space of splines
generated by the
$B_\alpha$'s, each vector $a_\alpha$ must lie in an affine subspace
$V_\alpha$, i.e. satisfy some
linear constraints. However there remain some degrees of freedom which are
used to minimize
 $\Vert a_\alpha \Vert_1$ for each $\alpha\in A$. It is easy to prove that
$\max \{\Vert a_\alpha \Vert_1 ; \alpha\in A\}$
is an upper bound of $\Vert Q_a \Vert_{\infty}$: thus, instead of
minimizing the infinite norm of $Q_a$,
 which is a difficult problem,
we minimize an upper bound of this norm, which is much easier to do.
Moreover, the latter problem has
always at least one solution, which is associated with a NB QI.
In the first part of the paper, we give a survey on NB univariate or
bivariate spline QIs
 defined on uniform or non-uniform partitions and already studied by the
author and coworkers.
In the second part, we give some new results, mainly on  univariate and
bivariate integral QIs on
{\sl non-uniform} partitions: in that case, NB QIs are more difficult to
characterize and the
optimal properties strongly depend on the geometry of the partition.
Therefore we have
restricted our study to QIs having interesting shape properties and/or
infinite norms uniformly
bounded independently of the partition.

%*********************************** INTRODUCTION   ***************************

\section{Introduction and notations}

The spline quasi-interpolants (abbr. QIs) considered in this paper have the
following general
form:
$$
Q_af=\sum_{\alpha\in A} \Lambda_\alpha (f) B_\alpha.
$$
$A$ denotes a finite or infinite set of indices. $ B_\alpha$ is a B-spline
with support $\Sigma_\alpha$,
defined on a uniform or nonuniform partition. The index $a$ of $Q_a$ refers
to a family of vectors
$\{a_\alpha, \alpha\in A\}$ which is described below.
The coefficients $\Lambda_\alpha (f)$ are discrete or integral functionals
of the following types:
$$
\Lambda_\alpha (f)=\sum_{\gamma \in \Gamma_\alpha} a_{\alpha} (\gamma)
f(x_\gamma) \;\; or \;\;
\Lambda_\alpha (f)=\sum_{\gamma \in \Gamma_\alpha} a_{\alpha} (\gamma)
\int_{\Sigma_\gamma}
B_\gamma^{*}f,
$$
where $\Gamma_\alpha$ is a finite set of indices nearby $\alpha$.
In the discrete coefficient functionals $\Lambda_\alpha (f)$, the points
$x_\gamma$ lie in
$\Sigma_\gamma$ and the associated QIs are called {\sl discrete
quasi-interpolants} (abbr.
dQIs). Those with integral coefficient functionals $\Lambda_\alpha (f)$ are
called {\sl integral
quasi-interpolants} (abbr. iQIs): the weight function in the integral is a
B-spline $B_{\gamma}^*$
(which can be different from $ B_\gamma$) satisfying $\int_{\Sigma_\gamma}
B_\gamma^{*}=1$ for all indices $\gamma$.
For $\alpha\in A$, the vectors $a_{\alpha}=(a_{\alpha}(\gamma),
\gamma \in \Gamma_\alpha)$  are determined in the following way:

(i) let $\PP$ be a space of polynomials included in the space of splines
generated by the family
$\{B_\alpha, \alpha\in A\}$, then we impose that $Q_a$ be exact on $\PP$,
i.e. $Q_a p=p$ for all
$p\in \PP$. This is equivalent to say that, for each $\alpha \in A$, the
vector $a_{\alpha}$ belongs
to some affine subspace $V_\alpha$ of
$\RR^{n_\alpha}$ where $n_\alpha=card(\Gamma_\alpha)$ is large enough.

(ii) in general there remain some undetermined components of $a_{\alpha}$.
Then we minimize the
$l^1$-norm of this vector: using standard results in optimization, it is
easy to prove that, for
each $\alpha\in A$, there exists at least one vector $a_{\alpha}\in
V_\alpha$ solution to this
problem.

The infinite norm of $Q_a$ is bounded above by $\nu(Q_a)=\max_{\alpha \in
A} \nu (a_\alpha)$
where $\nu (a_\alpha)=\Vert a_\alpha \Vert_1$. Thus, instead of minimizing
this norm, we
minimize $\nu(a_\alpha)$ for all $\alpha \in A$, the vector $a_\alpha$
satisfying the constraints
$a_\alpha\in V_\alpha$, which is much easier to do.  When $a_\alpha^*$ is a
solution of this
minimization problem, we say that the corresponding QI
$Q^*=Q_{a^*}$ is a {\sl near-best quasi-interpolant} (abbr. NBQI).
In the first part of the paper, we present various examples of NBQIs which
have been studied in
\cite{bis}\cite{mjib} and in the reports \cite{biss1}-\cite{biss3}. In the
second part,
we present  new results on some families of univariate and bivariate iQIs:
the proofs are somewhat
sketched and will be detailed elsewhere \cite{biss4}\cite{ms1}\cite{ms2}.

%************************ NB UNIVARIATE QI ON UNIFORM PARTITIONS
*************************

\section{Univariate Near-best QIs on uniform partitions}

Butzer et al. \cite{bssv}\cite{bs} have completed the results of Schoenberg
\cite{Sch1} about the
general expressions of spline dQIs of maximal approximation order (i.e.
exact on polynomials having
the same degree as the underlying spline) on the real line with $\ZZ$ as
sequence of knots. Another
technique, based on central difference operators, has been  given by the
author in \cite{Sab5}.

\subsection{Near-best spline dQIs}

Let $M_{2m}(x)$ denote the centered cardinal B-spline. Consider the family
of spline dQIs of order
$2m$ depending on $n+1$ arbitrary parameters $a=(a_0,a_1,\ldots,a_n)$, $n\ge m$:
$$
Q_af=\sum_{i\in\ZZ} \Lambda f(i) M_{2m}(x-i),
$$
with coefficient functionals
$$
\Lambda f(i)=a_0f(i)+\sum_{j=1}^{n}a_j\left(f(i+j)+f(i-j)\right).
$$
 Setting $\nu(a)=\vert a_0
\vert+\sum_{j=1}^{n}\vert a_j\vert$, then we have $\Vert
Q_a\Vert_{\infty}\le \nu(a)$.
By imposing that $Q_a$ be exact on the space $\Pi_r$ of polynomials of
degree at most $r$,
with $0\le r\le 2m-1$, we obtain a set of linear constraints: $a\in
V_r\subset {\RR}^{n+1}$.
We say that $Q^*=Q_{a^*}$ is a {\sl near best dQI} if
$$
\nu(a^*)=\min\{\nu(a); a\in V_r\}.
$$
There is existence, but in general not unicity, of solutions (see
\cite{biss1}, \cite{mjib}).

\noindent
{\sl Example:} cubic splines \cite{mjib}.
There is a unique optimal solution for $r=3$ and $n\ge 2$:
$$
a_0^*=1+\frac{1}{3n^2},\;\; a_n^*=-\frac{1}{6n^2},\;\; a_j^*=0 \;\; for
\;\; 1\le j\le n-1.
$$
Moreover, for all $n\ge 4$, $\Vert Q^* \Vert_{\infty}\le 1+\frac{2}{3n^2}.$
Here are the first values of $\Vert Q^* \Vert_{\infty}$ \& $\nu(a^*);$
$n=1: 1.222\;\; \&\;\; 1.666;$
$n=2: 1.139\;\;  \&\;\; 1.166;$
$n=3: 1.074\;\;  \&\;\; 1.074.$

%--------------------------- UNIVARIATE SPLINE iQIs ON UNIFORM PARTITIONS-----------------

\subsection{Near-best spline iQIs}

A similar study can be done for integral spline QIs. We refer to
\cite{biss1}\cite{mjib} and we only
give  an example given in these papers. Setting $a=(a_0,a_1,\ldots,a_n)$,
$n\ge m$,
$M_i(x)=M_{2m}(x-i)$ and $\langle f,M_i\rangle=\int fM_i$, we consider
$
Q_af=\sum_{i\in\ZZ} \Lambda f(i) M_i
$
with coefficient functionals
$$
\Lambda f(i)=a_0\langle f,M_i\rangle+\sum_{j=1}^{n}a_j\left(\langle
f,M_{i-j}\rangle+\langle
f,M_{i+j}\rangle\right).
$$
As in section 4.2, we have  $\Vert Q_a\Vert_{\infty}\le \nu(a)$
and we say that $Q^*=Q_{a^*}$ is a {\sl near best iQI} if
$
\nu(a^*)=\min\{\nu(a); a\in V_r\}.
$
There is existence, but in general not unicity, of solutions.

\noindent
{\sl Example:} cubic splines \cite{mjib}.
There is a unique optimal solution for $r=3$ and $n\ge 2$:
$$
a_0^*=1+\frac{2}{3n^2},\;\; a_n^*=-\frac{1}{3n^2},\;\; a_j^*=0 \;\; for
\;\; 1\le j\le n-1.
$$
Moreover, for all $n\ge 4$, $\Vert Q^* \Vert_{\infty}\le 1+\frac{4}{3n^2}.$
Here are the first values of $\Vert Q^* \Vert_{\infty}$ \& $\nu(a^*);$
$n=1: 1.5278\;\; \&\;\; 2.333;$
$n=2: 1.2778\;\;  \&\;\; 1.333;$
$n=3: 1.1481\;\;  \&\;\; 1.1482.$

%***************************  BIVARIATE NB QI ON UNIFORM PARTITIONS
*****************************

\section{Bivariate Near-best QIs on uniform partitions}

\subsection{A general construction of dQIs}

Let $\varphi$ be some bivariate B-spline on one of the two classical three or
four direction meshes of the plane (e.g. box-splines or H-splines, see
\cite{biss2}\cite{dB2}\cite{BHS}\cite{dBHR}\cite{Chu}\cite{Sab2}\cite{Sab5}\
cite{Sab6}\cite{Sb}).
Let
$\Sigma=supp(\varphi)$ and $\Sigma^*=\Sigma\cap{\ZZ}^2$.  Let $a$ be the
hexagonal (or
lozenge=rhombus) sequence formed by the values
$\{\varphi(i),i\in\Sigma^*\}$. The associated central difference operator
$\mathcal{D}$ is an isomorphism of $\PP(\varphi)$, the maximal subspace of
"complete " polynomials in the space of splines $\mathcal{S}(\varphi)$
generated by the integer translates of the B-spline $\varphi$ (see e.g.
\cite{BHS}\cite{dBHR}).
Computing the expansion of $a$ in some basis of the space of
hexagonal (or lozenge) sequences amounts to expand $\mathcal{D}$
in some basis of central difference operators.
Then, computing the formal inverse ${\mathcal{D}}^{-1}$ allows to define the dQI
$$
{\mathcal{Q}}f=\sum_{k\in {\ZZ}^2}{\mathcal{D}}^{-1}f(k)\varphi(\cdot-k)
$$
which is exact on ${\PP} (\varphi)$. Let us now give two examples of NB
dQIs which are detailed in
\cite{mjib}. The definition of these operators is quite the same as in
section 2.1.

\subsection{Near-best spline dQIs on a three direction mesh}

For example, let $\varphi$ be the $C^2$ quartic box-spline with support
$\Sigma=H_2$ where $H_s$
denotes the regular hexagon with edges of length $s\ge 1$, centered at the
origin  and
$H_s^*=H_s\cap {\ZZ}^2$. The near-best dQIs, which are exact on $\Pi_3$,
have coefficient
functionals with supports consisting of the center and of the 6 vertices of
$H_s^*, s\ge 1$. The
coefficients of values of
$f$ at those points are respectively
$1+\frac{1}{2s^2}$ and $-\frac{1}{12s^2}$, therefore the infinite norm of
the optimal dQIs
$Q_s^*$ is bounded above by $\nu_s^*=1+\frac{1}{s^2}$. Here are the first
values of
$\Vert Q^* \Vert_{\infty}$ \& $\nu_s^*;$
$n=1: 1.34028\;\; \&\;\; 2 ;$
$n=2: 1.22917\;\;  \&\;\; 1.25;$
$n=3: 1.10185\;\;  \&\;\; 1.111.$

\subsection{Near-best spline dQIs on a four direction mesh}

For example, let $\varphi$ be the $C^1$ quadratic box-spline. Let
$\Lambda_s$ be the lozenge
(rhombus) with edges of length $s\ge 1$, centered at the origin, and let
$\Lambda_s^*=\Lambda_s\cap {\ZZ}^2$. The near-best dQIs which are exact on
$\Pi_2$ have coefficient
functionals with supports consisting of the center and the 4 vertices of
$\Lambda_s^*, s\ge 1$. The
coefficients of values of $f$ at those points are respectively
$1+\frac{1}{2s^2}$ and $-\frac{1}{8s^2}$, therefore the infinite norm of
the optimal dQIs
$Q_s^*$ is bounded above by $\nu_s^*=1+\frac{1}{s^2}$. Here are the first values of
$\Vert Q^* \Vert_{\infty}$ \& $\nu_s^*;$
$n=1: 1.5\;\; \&\;\; 2 ;$
$n=2: 1.25\;\;  \&\;\; 1.25;$
$n=3: 1.111\;\;  \&\;\; 1.111.$

Examples with $C^2$ quartic box-splines are also given in \cite{mjib}.

%**************************  UNIVARIATE NB QI ON NON UNIFORM PARTITIONS
****************

\section{Univariate Near-best QIs on non-uniform partitions}

In this section, we consider dQIs or iQIs of degree $m\ge 2$ defined on a
bounded interval
$I=[a,b]$. Let $E_m=\{-m+1,\ldots,0\}$ and $J=\{0,1,\ldots, m+n-1\}$. Let
$T=\{t_i; i\in E_m\cup
J\}$ be an arbitrary non-uniform increasing sequence of knots
with multiple knots at $t_0=a$ and $t_n=b$, as usual. Let $B_j$ be the
B-spline with support
$[t_{j-m},t_{j+1}]$ for $j\in J$, and let $e_p(x)=x^p$ for all $p\ge 0$.
Setting $\theta_j=\frac1m \sum_{s\in E_m}t_{j+s}$, it is well known that
$e_1=\sum_{j\in J}\theta_j B_j$ and
$e_2=\sum_{j\in J}\theta_j^{(2)} B_j$, with
$\theta_j^{(2)}=\frac{2}{m(m-1)}\sum_{(r,s)\in E_m^2,r<s}t_{j+r}t_{j+s}$.
We recall the following expansion \cite{ms}:
$$
\lambda_j=\theta_j^2-\theta_j^{(2)}=\frac{1}{m^2 (m-1)}\sum_{(r,s)\in E_m^2,r<s}
(t_{j+r}-t_{j+s})^2>0.
$$
More generally $e_r=\sum_{j\in J}\theta_j^{(r)} B_j$ where the
$\theta_j^{(r)}$'s are proportional
to symmetric functions of knots.
The simplest dQI is the Schoenberg operator
$$
S_1f=\sum_{i\in J} f(\theta_i) B_i
$$
which is exact on $\Pi_1$ and is also shape preserving. Moreover, it satisfies
$$
S_1e_2-e_2=\sum_{i\in J} \lambda_j B_i\ge 0.
$$

\subsection{Uniformly bounded dQIs}

Let us only give an example: we start from a family of differential QIs of
degree $m$ which are {\sl
exact on} $\Pi_2$ (see also \cite{biss3}):
$$
S_2^*f=\sum_{j\in J} \lambda_j^{(2)}(f) B_j,\quad
\lambda_j^{(2)}(f)=f(\theta_j)-\frac12 \lambda_j D^2f(\theta_j).
$$

On the other hand, $\frac12 D^2f(\theta_j)$ can be replaced on the space
$\Pi_2$ by the
second order divided difference $[\theta_{j-1}, \theta_j,\theta_{j+1}]f$,
therefore the dQI
defined by
$$
S_2f=\sum_{j\in J} \mu_j^{(2)}(f) B_j,\;\;
\mu_j^{(2)}(f)=f(\theta_j)- \lambda_j[\theta_{j-1}, \theta_j,\theta_{j+1}]f
$$
is also {\sl exact on} $\Pi_2$. Moreover, one can write
$$
\mu_i^{(2)}(f)=a_if_{i-1}+b_i f_i+c_i f_{i+1}
$$
with
$
a_i=-\lambda_i/\Delta \theta_{i-1}(\Delta \theta_{i-1}+\Delta \theta_{i}),
\;\;
c_i=-\lambda_i/\Delta \theta_{i}(\Delta \theta_{i-1}+\Delta \theta_{i}),$
 and \\
$b_i=1+ \lambda_i/\Delta \theta_{i-1}\Delta \theta_{i}$.
So, according to the introduction
$$
\Vert S_2 \Vert_{\infty}\le \max_{i\in J} (\vert a_i \vert+\vert b_i \vert
+\vert c_i \vert)
\le 1+2\max_{i\in J}\frac{\lambda_i}{\Delta \theta_{i-1}\Delta \theta_{i}}.
$$
The following theorem \cite{biss3} extends a result given for quadratic
splines in
\cite{mjib}\cite{Sab4}.\\

\noindent
\begin{theorem}
For any degree $m$, the dQIs $S_2$ are uniformly bounded independently of
the partition. More
specifically, if $[r]$ denotes the floor of $r$:
$$
\Vert S_2 \Vert_{\infty}\le [\frac12 (m+4)]
$$
\end{theorem}

{\bf Remark.} For quadratic splines, one can prove that $\Vert S_2
\Vert_{\infty}\le 2.5$ for all
partitions. For uniform partitions, one gets $\Vert S_2
\Vert_{\infty}=\frac{305}{207}\approx
1.4734$

%---------------------------------------------------------------------------

\subsection{Near-best dQIs}

Let us consider the family of dQIs of degree $m$ defined by
$$
Qf=Q_{p,q}f=\sum_{i\in\ZZ} \mu_i(f) B_i.
$$
depending on the two integer parameters $p\ge m$ and $q\le \min(m,2p)$.
Their coefficient functionals depend on $2p+1$ parameters
$$
\mu_i(f)=\sum_{s=-p}^{p}\lambda_i(s)f(\theta_{i+s}),
$$
and we impose that $Q$ is {\sl exact on the space} $\Pi_q$.
The latter condition is equivalent to $Q e_r=e_r$ for all monomials of
degrees $0\le r\le q$.
It implies that for all indices $i$, the parameters $\lambda_i(s)$ satisfy
the system of $q+1$
linear equations:
$$
\sum_{s=-p}^{p}\lambda_i(s)\theta_{i+s}^r=\theta_i^{(r)},\quad 0\le r\le q.
$$
The matrix $V_i\in {\RR}^{(q+1)\times (2p+1)}$ of this system, with
coefficients
$V_i(r,s)=\theta_{i+s}^r$, is a Vandermonde
matrix of maximal rank $q+1$, therefore there are $2p-q$ {\sl free parameters}.
Denoting by $b_i\in {\RR}^{q+1}$ the vector in the right hand side, with
components
$b_i(r)=\theta_i^{(r)},\;\; 0\le r\le q$, we consider the sequence of
minimization problems,
for $i\in \ZZ$:
$$
\min \{\Vert \lambda_i \Vert_1;\;\; V_i \lambda_i=b_i\}.
$$
We have seen in the introduction that $\nu_1^*(Q)=\max_{i\in\ZZ}\min\Vert
\lambda_i \Vert_1$ is an upper bound
of $\Vert Q_q \Vert_{\infty}$ which is easier to evaluate than the true
norm of the dQI.
The {\sl objective function being convex} and {\sl the domains being affine
subspaces}, we have the
following

\begin{theorem} The above minimization problems have always solutions,
which, in general,
are non unique.
\end{theorem}

\noindent
Let us give an example of optimal dQIs (\cite{bis}\cite{biss3}\cite{mjib})
in order to bring out
the constraints on the partition implied by the optimal property. For $m=2$
and $p\ge 2$, assume
that the partition $T$ satisfies, for all $i$
$$
\theta_{i-1}+\theta_i\le \theta_{i-p}+\theta_{i+p}\le \theta_i+\theta_{i+1},
$$
then there is a unique optimal solution (here $h_i=t_i-t_{i-1}$):
$$
\lambda_i^*(-p)=-\frac14\frac{h_i^2}{(\theta_{i-p}-\theta_{i+p})(\theta_{i}-
\theta_{i+p})},\;\;
\lambda_i^*(p)=-\frac14\frac{h_i^2}{(\theta_{i-p}-\theta_{i+p})(\theta_{i+p}
-\theta_{i})}
$$
$$
\lambda_i^*(p)=1+\frac14\frac{h_i^2}{(\theta_{i+p}-\theta_{i})(\theta_{i}-
\theta_{i-p})},
\;\;  \lambda_i^*(s)=0 \;\; for \;\; s\neq 0,-p,p.
$$
Denote by $Q_{p,2}^*$ the associated dQI, then we have \cite{biss3}

\begin{theorem} The infinite norm of $Q_{p,2}^*$ is uniformly bounded
independently of $p$ and
of the partition $T$:
$$
\Vert Q_{p,2}^* \Vert_{\infty} \le 3.
$$
\end{theorem}

%---------------------------------------------------------------------------
---------------------
\subsection{Uniformly bounded iQIs}

Various types of integral QIs  are considered in
\cite{Cz}\cite{SSb}\cite{biss4}.
Here, we restrict our study  to Goodman-Sharma (abbr. GS)
type iQIs which  appear in \cite{GS1}. They are simpler than those we have
already
studied in \cite{SSb}. Given $I=[a,b]$ and the sequence of knots
$T=\{t_{-m}=\ldots=t_0=a<t_1<\ldots<t_{n-1}<b=t_n=\ldots=t_{n+m}\}$, let
$B_i$ be the B-spline of
degree $m$ and support $[t_{i-m},t_{i+1}]$ normalized by
$\sum_{i=0}^{n+m-1}B_i=1$, and let $\tilde
M_{i-1}(t)$ be the B-spline of degree $m-2$ with support
$\tilde \Sigma_{i-1}=[t_{i-m+1},t_i]$, normalized by $\tilde
\mu_{i}^{(0)}=\tilde \mu_{i}(e_0)=
\int_a^b\tilde M_{i-1}(t)=1$.
The original GS-type iQI can be written as follows
$$
G_1f=\sum_{i=0}^{n+m-1}\tilde\mu_i(f)B_i,
$$
where $\tilde\mu_0(f)=f(t_0)$, $\tilde\mu_{n+m-1}(f)=f(t_n)$ and, for $1\le
i\le n+m-2$,
$$
\tilde\mu_i(f)=\int_a^b \tilde M_{i-1}(t)f(t)dt.
$$
 It is easy to verify that $G_1$ is {\sl exact on} $\Pi_1$ and that $\Vert
G_1\Vert_{\infty}=1$. It has also interesting {\sl shape preserving
properties}: for example it is
proved in \cite{GS1} that  $G_1$  preserves the {\sl positivity and the
convexity} of the
approximated function $f$. Moreover, the authors also prove that
$$
G_1e_2-e_2=\frac{2m}{m+1}(S_1 e_2-e_2).
$$
Let us also prove that it preserves the {\sl monotonicity} of $f$.

\begin{theorem} If $f$ is a monotone function on $I$, then $G_1f$ is also
monotone with the same
sense of variation.
\end{theorem}
{\sl Proof.} Let us use the simplified notation
$G_1f=\sum_{i=0}^{n+m-1}\tilde\mu_i B_i$. Then the
first derivative is
$DG_1f=m\sum_{i=1}^{n+m-1}\frac{\tilde\mu_i-\tilde\mu_{i-1}}{t_i-t_{i-m}}
B_i^*$, where $B_i^*$ is the B-spline of degree $m-1$, normalized as $B_i$.
Assuming that $f$ is
monotone increasing, we have to prove that $DG_1f\ge 0$. For that, it is
sufficient to prove that
the sequence $\{(\tilde\mu_i-\tilde\mu_{i-1}), 1\le i\le n+m-1\}$ is
monotone increasing. For
$i=1$, we have $\tilde\mu_1-\tilde\mu_0=\int_{t_0}^{t_1}\tilde M_0 f-f(t_0)$; as
$\int_{t_0}^{t_1}\tilde M_0=1$, the mean value theorem gives
$\int_{t_0}^{t_1}\tilde M_0
f=f(\xi_0)$ for some $\xi_0\in [t_0,t_1]$, hence $\tilde\mu_1-\tilde\mu_0=f(\xi_0)-f(t_0)\ge 0$.
The same kind of proof holds for
$i=n+m-1$. Now, for $2\le i\le n+m-2$, we know that $\tilde M_i-\tilde
M_{i-1}=-DB_{i-2}^*=$
(see e.g. \cite{dB1}), therefore
$\tilde\mu_i-\tilde\mu_{i-1}=\frac{1}{t_i-t_{i-m}}\int
(-DB_{i-2}^*f)$ and, after integration by parts,
$\tilde\mu_i-\tilde\mu_{i-1}=\frac{1}{t_i-t_{i-m}}\int B_{i-2}^*f'\ge 0.$\\

Therefore, we can see that the operator $G_1$ is very close to the Schoenberg's
operator.
>From the expression of $G_1e_2-e_2$, we can deduce the family of GS-type
iQIs defined by
$$
G_2f=f(t_0)B_0+\sum_{i=1}^{n+m-2}[a_i\tilde\mu_{i-1}(f)+b_i\tilde\mu_i(f)+c_
i\tilde\mu_{i+1}(f)]B_i
+f(t_n)B_{n+m-1},
$$
which are {\sl exact} on $\Pi_2$. The three constraints
$G_2e_k=e_k,\;\;k=0,1,2$, lead to
the  following system of equations, for $1\le i\le n+m-2$:
$$
a_i+b_i+c_i=1,\quad
\theta_{i-1}a_i+\theta_{i}b_i+\theta_{i+1}c_i=\theta_{i},\quad
\tilde\mu_{i-1}^{(2)}a_i+\tilde\mu_{i}^{(2)}b_i+\tilde\mu_{i+1}^{(2)}c_i=
\theta_i^{(2)}.
$$
This is a consequence of the following facts
$$
\tilde\mu_i(e_1)=\int_a^b t\tilde
M_{i-1}(t)dt=\frac{1}{m}\sum_{s=1}^{m}t_{i-m+s}=\theta_i,
$$
$$
\tilde\mu_i^{(2)}=\mu_i(e_2)=\int_a^b t^2\tilde
M_{i-1}(t)dt=\frac{2}{m(m+1)}\tilde s_2(T_i)
$$
$$
=\frac{2}{m(m+1)}\sum_{1\le r\le s\le m}t_{i-m+r}t_{i-m+s}
$$
For $m=2$ (quadratic splines), here are the explicit expressions of the
coefficients:
$$
a_i=\frac{-h_i^2}{(h_{i-1}+h_i)(h_{i-1}+h_i+h_{i+1})},\;\;
c_i=\frac{-h_i^2}{(h_{i-1}+h_i+h_{i+1})(h_i+h_{i+1})},
$$
and $b_i=1-a_i-c_i$. It is clear that
$\Vert G_2 \Vert_{\infty}\le 1+2\max (\vert a_i\vert +\vert c_i\vert)\le 5$.
In fact, this upper bound is valid for any degree $m$: the computation of
coefficients and the
proof of this interesting result will be given in \cite{biss4}.

\begin{theorem}
The iQIs $G_2$ are uniformly bounded independently of the partition and of
the degree $m$. More specifically, one has
$$
\Vert G_2 \Vert_{\infty}\le 5
$$
\end{theorem}

%********************  BIVARIATE QUADRATIC QI ON CRISS-CROSS TRIANGULATIONS
*********************

\section{Bivariate quadratic spline QIs on non-uniform partitions}

At the author's knowledge, the only bivariate box-splines which have been
extended to non uniform
partitions of the plane are $C^1$-quadratic box-splines on criss-cross
triangulations
\cite{CW}\cite{Sab1}. Recently \cite{Sab7}-\cite{Sab10} we have constructed
a set of B-splines
generating the space of quadratic splines on a rectangular domain and {\sl
having their support in
this domain}. Moreover, we have
defined a discrete quasi-interpolant which is exact on $\Pi_2$ and
uniformly bounded
independently of the partition (it is different from the operator
introduced in \cite{CW}, see e.g.
\cite{DaL},\cite{DaS}). For the sake of simplicity, we assume here that the
domain is the whole
plane endowed with a nonuniform criss-cross triangulation obtained by
drawing diagonals in each
rectangle
$R_{ij}=[x_{i-1},x_i]\times[y_{j-1},y_j]$ of the partition defined by the
two sequences $X=\{x_i,\;
i\in \ZZ\}$ and  $Y=\{y_j,\;j\in \ZZ\}$. We set
$h_i=x_i-x_{i-1}$,
$k_j=y_j-y_{j-1}$, $s_i=\frac12(x_{i-1}+x_i)$ and
$t_j=\frac12(y_{j-1}+y_j)$. Let $B_{ij}$ be the
B-spline whose octagonal support $\Sigma_{ij}$ is centered at the point
$\omega_{ij}=(s_i,t_j)$.
Let $\Pi_{ij}$ be the continuous piecewise affine pyramid (of egyptian
type) satisfying
$\Pi_{ij}(\omega_{ij})=1$, whose support is the central rectangle $R_{ij}$
of $\Sigma_{ij}$. As
$\int_{R_{ij}}\Pi_{ij}=\frac13h_ik_j$, we can also define the pyramid $\tilde
\Pi_{ij}=\frac{3}{h_ik_j}\Pi_{ij}$ normalized by
$\int_{R_{ij}}\tilde\Pi_{ij}=1$.
Finally, in the same way, let $\chi_{ij}$ be the characteristic function of
$R_{ij}$ and let
$\tilde \chi_{ij}=\frac{1}{h_ik_j}\chi_{ij}$ normalized by
$\int_{R_{ij}}\tilde\chi_{ij}=1$.
For monomials, we use the notation $e_{rs}(x,y)=x^r y^s$.
We know study some families of discrete and integral quasi-interpolants.
Detailes proofs will be
given elsewhere.

%------------------------------------------------------------------------------------------------
\subsection{Discrete quasi-interpolants}

The simplest dQI is the Schoenberg operator defined by
$$
S_1f=\sum_{(i,j)\in\ZZ^2}f(\omega_{ij})B_{ij}.
$$
It is well known (\cite{CW}, \cite{DaL}) that  $S_1e_{rs}=e_{rs}$ for $0\le
r,s \le 1$ and
$$
S_1e_{20}=e_{20}+\frac14\sum_{(i,j)\in\ZZ^2} h_i^2 B_{ij},\;\;
S_1e_{02}=e_{02}+\frac14\sum_{(i,j)\in\ZZ^2} k_j^2 B_{ij}.
$$
$S_1$ is clearly a positive operator and it preserves the bimonotonicity
and the biconvexity of
$f$ (i.e. the monotonicity and the convexity in the directions of
coordinate axes). This can be
proved by using the expressions of partial derivatives $\partial_1 B_{ij}$
and $\partial_2
B_{ij}$ which are piecewise affine functions whose values are given in the
technical report
\cite{Sab10}.

%---------------------------------------------------------------------------
---------------------
\subsection{Integral quasi-interpolants}

We now study the two following iQIs:
$$
T_1f=\sum_{(i,j)\in\ZZ^2}\langle f, \tilde\Pi_{ij}\rangle B_{ij}
$$
where $\langle f,g \rangle =\int_{\RR^2} f(v)g(v) dv$, and the GS-type iQI:
$$
G_1f=\sum_{(i,j)\in\ZZ^2}\langle f, \tilde\chi_{ij}\rangle B_{ij}.
$$
Let $\mu_{ij}(f)=\frac{1}{h_ik_j}\int_{R_{ij}}f$ be the mean value of $f$
on $R_{ij}$, then one
can  also  write
$$
G_1f=\sum_{(i,j)\in\ZZ^2}\mu_{ij}(f) B_{ij}.
$$
These two operators are very close to each other and to the Schoenberg
operator defined before.
In particular, they are  exact on bilinear polynomials and they satisfy
respectively
$$
T_1e_{20}=e_{20}+\frac{3}{10}\sum_{(i,j)\in\ZZ^2} h_i^2 B_{ij},\;\;
T_1e_{02}=e_{02}+\frac{3}{10}\sum_{(i,j)\in\ZZ^2} k_j^2 B_{ij},
$$
$$
G_1e_{20}=e_{20}+\frac13\sum_{(i,j)\in\ZZ^2} h_i^2 B_{ij},\;\;
G_1e_{02}=e_{02}+\frac13\sum_{(i,j)\in\ZZ^2} k_j^2 B_{ij}.
$$
Moreover, as $S_1$, they both preserve the bimonotonicity and the
biconvexity of $f$.

>From the properties of $T_1$ and $G_1$ on monomials $e_{20}$ and $e_{02}$,
one can deduce the
two following iQIs which are both exact on the space $\Pi_2$ of bivariate
quadratic polynomials:
$$
T_2f=\sum_{(i,j)\in\ZZ^2}\langle f, M_{ij}\rangle B_{ij},\;\;
G_2f=\sum_{(i,j)\in\ZZ^2}\langle f, \psi_{ij}\rangle B_{ij}.
$$
where the two functions $M_{ij}$ and $\psi_{ij}$ are respectively defined by
$$
M_{ij}=a_i\tilde\Pi_{i-1,j}+\bar
a_i\tilde\Pi_{i+1,j}+b_{ij}\tilde\Pi_{ij}+c_j\tilde\Pi_{i,j-1}+\bar
c_j\tilde\Pi_{i,j+1}
$$
$$
\psi_{ij}=\alpha_i\tilde\chi_{i-1,j}+\bar
\alpha_i\tilde\chi_{i+1,j}+\beta_{ij}\tilde\chi_{ij}+\gamma_j\tilde\chi_{i,j
-1}+\bar
\gamma_j\tilde\chi_{i,j+1}
$$
The coefficients of these functions are the following
$$
a_i=\frac{-3h_i^2}{(h_{i-1}+h_i)(3h_{i-1}+4h_i+3h_{i+1})},\;
\bar a_i=\frac{-3h_i^2}{(3h_{i-1}+4h_i+3h_{i+1})(h_i+h_{i+1})}
$$
$$
c_j=\frac{-3k_j^2}{(k_{j-1}+k_j)(3k_{j-1}+4k_j+3k_{j+1})},\;
\bar c_j=\frac{-3k_j^2}{(3k_{j-1}+4k_j+3k_{j+1})(k_j+k_{j+1})},
$$
and $b_{ij}=1-(a_i+\bar a_i+c_j+\bar c_j)$, for the operator $T_2$.
It is easy to verify that for all criss-cross triangulations, one has
$
\vert a_i \vert, \vert \bar a_i \vert, \vert c_j \vert, \vert \bar c_j
\vert \le \frac34
$.
For the operator $G_2$, we get
$$
\alpha_i=\frac{-h_i^2}{(h_{i-1}+h_i)(h_{i-1}+h_i+h_{i+1})},\;
\bar \alpha_i=\frac{-h_i^2}{(h_{i-1}+h_i+h_{i+1})(h_i+h_{i+1})}
$$
$$
\gamma_j=\frac{-k_j^2}{(k_{j-1}+k_j)(k_{j-1}+k_j+k_{j+1})},\;
\bar \gamma_j=\frac{-k_j^2}{(k_{j-1}+k_j+k_{j+1})(k_j+k_{j+1})},
$$
and $\beta_{ij}=1-(\alpha_i+\bar \alpha_i+\gamma_j+\bar \gamma_j)$.
It is easy to verify that for all criss-cross triangulations, one has
$
\vert \alpha_i \vert, \vert \bar \alpha_i \vert, \vert \gamma_j\vert, \vert
\bar \gamma_j \vert
\le 1.
$
>From these inequalities and the fact that
$$
\Vert T_2 \Vert_{\infty}\le 1+2\max\{\vert a_i \vert,\vert\bar
a_i\vert,\vert c_j\vert,\vert\bar
c_j \vert\},\;\;
\Vert G_2 \Vert_{\infty}\le 1+2\max\{\vert \alpha_i\vert,\vert \bar
\alpha_i\vert,\vert \gamma_j\vert,\vert \bar\gamma_j\vert\},
$$
one can deduce the following interesting result.

\begin{theorem} For all non-uniform criss-cross triangulations, the operators $T_2$
and $G_2$ are uniformly bounded. Moreover,
$$
\Vert T_2 \Vert_{\infty}\le 7 \;\; and \;\; \Vert G_2 \Vert_{\infty}\le 9.
$$
\end{theorem}

{\bf Remark}. In the case of {\sl uniform} criss-cross triangulations, one
gets respectively:
$$
a_i=\bar a_i=c_j=\bar c_j=\frac{-3}{20},\; b_{ij}=\frac85, \;
\Longrightarrow \;
\Vert T_2 \Vert_{\infty}\le \frac{11}{5}\approx 2.22
$$
$$
\alpha_i=\bar \alpha_i=\gamma_j=\bar\gamma_j=\frac{-1}{6},\;
\beta_{ij}=\frac53,
\;\Longrightarrow \;
\Vert G_2 \Vert_{\infty}\le \frac73\approx 2.33.
$$

\subsection{Powell-Sabin quasi-interpolants}

Recently (\cite{ms1},\cite{ms2}), using an interesting result by Dierckx
\cite{Dx}, we have
introduced and studied new families of quadratic splines quasi-interpolants
defined on Powell-Sabin
type triangulations. As for the previous QIs of this section, we have
obtained some families of QIs
which are both {\sl exact on} $\Pi_2$ and {\sl uniformly bounded}
independently of the partition.

%***************************************************************************
*******************

\section{Some applications}.

\subsection{Approximation of functions}

>From a classical result in approximation theory (see e.g. \cite{DVL}) we
know that if $Q$ is an
operator defined on a space of smooth functions $f$ with values in a space
of splines $\mathcal{S}$,
one has
$$
\Vert f-Qf \Vert_{\infty}\le (1+\Vert Q \Vert_{\infty})
d_{\infty}(f,\mathcal{S})
$$
As the various QIs studied above have uniformly bounded norms, their
approximation order
is only governed by $d_{\infty}(f,\mathcal{S})$, i.e. by the distance of
$f$ to the maximal space of
polynomials included in $\mathcal{S}$. Since the values of $\Vert Q
\Vert_{\infty}$ are  small, we
obtain quite good approximants which can be used in various fields of
numerical anlysis.

\subsection{Approximation of zeros}

Quadratic dQIs are simple and good approximants and their zeros are rather
easy to compute, so they give good approximations of the zeros of the
approximated function. We
already did some computations with orthogonal polynomials, and the first
results are encouraging
(see e.g.\cite{Sab7}, \cite{Sab9}).

\subsection{Quadrature formulas}

For the same reasons, quadrature formulas (QF) are easily obtained by
integrating spline QIs.
An interesting  univariate example is given in \cite{Sab7}. The study of
bivariate and trivariate QF
is still in progress and the results already obtained in two and three
variables are also
encouraging \cite{DeS}.

\subsection{Pseudo-spectral methods associated with quasi-interpolants}

Derivatives of QIs give quite good approximations of derivatives of the
approximated function. This
simple fact is the basic idea for developing pseudo-spectral methods based
on univariat and
multivariate dQIs with low degrees.

\medskip

%************************************  REFERENCES *************************

\textbf{Keywords:} spline operators, quasi-interpolants.

\textbf{AMS Classification:} 41A36.

%***************************************************************************
************

\bigskip

\footnotesize{P. Sablonni\`ere,
INSA, 20 avenue des Buttes de Co\"esmes,
CS 14315, F-35043-Rennes c\'edex, France.

email: psablonn@insa-rennes.fr}

\end{document}